\documentclass[reqno]{amsart}

\usepackage{amsfonts,amssymb,graphicx,float}

\restylefloat{figure}

\makeatletter                                                  %
\def\th@plain{\slshape}                                        %
\makeatother                                                   %

\newcommand{\oi}{[0,1]}
\newcommand{\Nbb}{\mathbb{N}}
\newcommand{\Zbb}{\mathbb{Z}}
\newcommand{\Qbb}{\mathbb{Q}}
\newcommand{\Rbb}{\mathbb{R}}

\newcommand{\ve}{\varepsilon}

\newcommand{\Dcal}{\mathcal{D}}

\newcommand{\newword}[1]{\textsl{#1}}
\newcommand{\vect}[3]{#1_#2,\ldots ,#1_#3}
\newcommand{\abs}[1]{\lvert#1\rvert}

\newcommand{\vc}[2]{(#1\thickspace #2)^{t}}

\DeclareMathSymbol{\upharpoonright}{\mathrel}{AMSa}{"16}

\DeclareMathSymbol{\nmid}{\mathrel}{AMSb}{"2D}

\DeclareMathOperator{\length}{length}
\DeclareMathOperator{\GL}{GL}

\theoremstyle{plain}
\newtheorem{theorem}{Theorem}[section]

\newtheorem{observation}[theorem]{Observation}

\theoremstyle{definition}

\newtheorem{example}[theorem]{Example}

\begin{document}

\bibliographystyle{plain}

\sloppy

\title[Lagrange Theorem]{A general Lagrange Theorem}

\author[G. Panti]{Giovanni Panti}
\address{Department of Mathematics\\
University of Udine\\
via delle Scienze 208\\
33100 Udine, Italy}
\email{panti@dimi.uniud.it}


\thanks{\emph{2000 Math.~Subj.~Class.}: 11J70}

\maketitle

\section{Introduction}

The ordinary continued fractions expansion of a real number is based on the Euclidean division. Variants of the latter yield variants of the former, all encompassed by a more general Dynamical Systems framework. For all these variants the Lagrange Theorem holds: a number has an eventually periodic expansion if and only if it is a quadratic irrational. This fact is surely known for specific expansions, but the only proof for the general case that I could trace in the literature follows as an implicit corollary from much deeper results by Boshernitzan and Carroll on interval exchange transformations~\cite{boshernitzancarroll97}. It may then be useful to have at hand a simple and virtually computation-free proof of a general Lagrange Theorem.

Let $Q\subset\Qbb\cap\oi$ be a finite or countable set containing $0$,
$1$, and at least another rational, and assume that every interval having as extrema two consecutive elements $p/q<r/s$ of $Q$ is
\newword{unimodular} (i.e., 
$\bigl\lvert \begin{smallmatrix}
p&r\\ q&s
\end{smallmatrix} \bigr\rvert=-1$). 
I always write rational numbers in reduced form; the unimodularity condition amounts then to saying that the column vectors $\vc pq$ and $\vc rs$ constitute a $\Zbb$-basis for $\Zbb^2$.
Let $\Dcal$ be the family of all half-open intervals $\Delta_a=(p/q,r/s]$, for $p/q,r/s\in Q$ as above; we agree that the index $a$ varies in a fixed set $I\subseteq\Zbb$.
We also agree that, if $Q\setminus\{1\}$ has maximum element $p/q$, then the open interval $(p/q,1)$ replaces $(p/q,1]$ in $\Dcal$; we then have $\bigcup\Dcal=(0,1)$. Fix an arbitrary function $\ve:I\to\{-1,+1\}$, and consider the matrix
$$
G_a=
\begin{pmatrix}
0 & 1 \\ 1 & 1
\end{pmatrix}
\begin{pmatrix}
0 & 1 \\ 1 & 0
\end{pmatrix}^{\frac{\ve(a)+1}{2}}
\begin{pmatrix}
p & r \\ q & s
\end{pmatrix}^{-1}.
$$

Writing elements $x\in(0,1)$ as vectors in projective coordinates $\vc x1$, multiplication to the left by $G_a$ induces a fractional-linear homeomorphism ---by abuse of language still denoted by $G_a$--- from $\Delta_a$ to either $(0,1]$, or $[0,1)$, or $(0,1)$;
$G_a$ is monotonically increasing or decreasing according whether $\ve(a)$ equals $-1$ or $+1$. The resulting piecewise-fractional map $G:(0,1)\to\oi$, defined by $G(x)=G_a(x)$ for $x\in\Delta_a$, is the \newword{Gauss map} determined by $Q$ and $\ve$.

Denote by $\psi_a$ the $a$th inverse branch of $G$; the domain of $\psi_a$ is $\oi$, and its range is the topological closure of $\Delta_a$.
Computing the inverse of $G_a$, we obtain explicitly
$$
\psi_a(x)=
\begin{cases}
\displaystyle{\frac{(r-p)x+p}{(s-q)x+q},}
&\text{if $\ve(a)=-1$;}\\
\displaystyle{\frac{(p-r)x+r}{(q-s)x+s},}
&\text{if $\ve(a)=+1$.}
\end{cases}
$$
If $x,Gx,G^2x,\ldots,G^{n-1}x$ are all in $(0,1)$ then, letting $G^{t-1}x\in\Delta_{a_t}$ for $1\le t\le n$, we have the identity
\begin{equation}\label{eq1}
x=\psi_{a_1}\psi_{a_2}\cdots\psi_{a_n}(G^nx).
\end{equation}
In Dynamical Systems language, the (finite or infinite) sequence $a_1,a_2,\ldots$ is the \newword{kneading sequence} of $x$. By~(\ref{eq1}), if $x$ has a finite kneading sequence (i.e., $G^nx\in\{0,1\}$ for some $n$), then $x$ is rational. Conversely, let $x=u/v\in\Delta_a$ be rational, and let $p/q<r/s$ be the extrema of $\Delta_a$. Then $\vc{u}{v}=l\vc{p}{q}+m\vc{r}{s}$ for some nonnegative integers $l,m$ with $m\ge1$. Multiplying by $G_a$ to the left, we see that the denominator of $Gx$ equals $l+m$, which is strictly less than the denominator $v=lq+ms$ of $x$. Hence the sequence of the denominators along the $G$-orbit of $x$ must end up at $1$, and $x$ must end up either at $0$ or at $1$.

\begin{example}
The first and main example is of course given by the ordinary continued fractions: $I=\{1,2,\ldots\}$, $\Delta_a=\bigl(1/(a+1),1/a\bigr]$, and $\ve=+1$ throughout. The Gauss map is $G(x)=1/x-\lfloor 1/x\rfloor$, and has graph
\begin{figure}[H]
\begin{center}
\includegraphics[height=3cm,width=3cm]{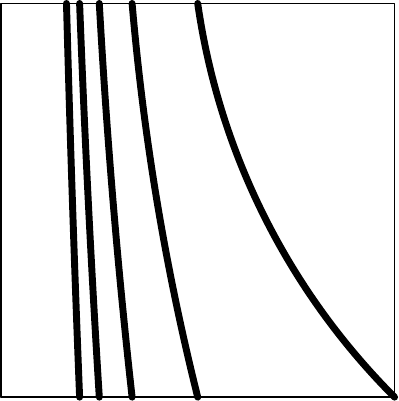}
\end{center}
\end{figure}
\noindent We have $\psi_a(x)=1/(a+x)$, and hence (1) assumes the familiar shape
$$
x=
\cfrac{1}{a_1+ 
\cfrac{1}{a_2+ 
\cfrac{1}{
\cfrac{\ddots}{a_{n-1}+
\cfrac{1}{a_n+G^nx} 
}}}}
$$
\end{example}

\begin{example}
In the odd continued fractions we have $I$ and the $\Delta_a$'s as above, $\ve(a)=(-1)^{a+1}$. The Gauss map is $G(x)=\abs{1/x-b}$, where $b$ is the odd integer in $\{\lfloor 1/x\rfloor,\lfloor 1/x\rfloor+1\}$, and has graph
\begin{figure}[H]
\begin{center}
\includegraphics[height=3cm,width=3cm]{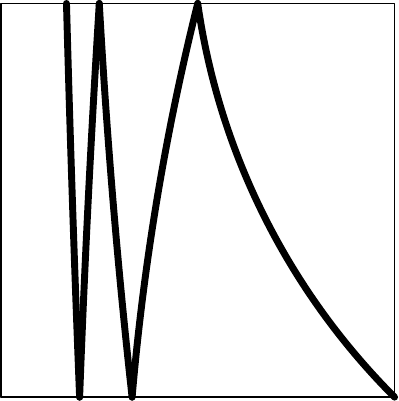}
\end{center}
\end{figure}
\noindent Setting $(b_t,\ve_t)$ equal to $(a_t,+1)$ for $a_t=\lfloor1/G^{t-1}x\rfloor$ odd, and to $(a_t+1,-1)$ for $a_t$ even, we have $\psi_{a_t}(x)=1/(b_t+\ve_tx)$, and therefore
$$
x=
\cfrac{1}{b_1+ 
\cfrac{\ve_1}{b_2+ 
\cfrac{\ve_2}{
\cfrac{\ddots}{b_{n-1}+
\cfrac{\ve_{n-1}}{b_n+\ve_nG^nx} 
}}}}
$$
\end{example}

Many more examples can be found in~\cite{baladivallee05}, \cite{vallee06}, and references therein, along with a detailed analysis of the ergodic properties of these Gauss maps.

\section{Convergence $\ldots$}

Unimodular intervals have a curious property.

\begin{observation}
Let\label{ref1} $\Gamma_1\supset\Gamma_2\supset\cdots$ be a nested sequence of closed unimodular intervals, each one strictly contained in the previous one. Then
$$
\lim_{n\to\infty}\length(\Gamma_n)=0.
$$
\end{observation}
\begin{proof}
The length of $\Gamma_n=[p_n/q_n,r_n/s_n]$ is $(q_ns_n)^{-1}$, which is less than or equal to $[\max(q_n,s_n)]^{-1}$. If $u/v$ is a rational in the topological interior of $\Gamma_n$, then $v$ is strictly greater than $\max(q_n,s_n)$ (again, because $\vc{u}{v}=l\vc{p_n}{q_n}+m\vc{r_n}{s_n}$, for certain $0<l,m\in\Zbb$). It follows that the sequence $\max(q_1,s_1), \max(q_2,s_2),\ldots$ is strictly increasing and goes to infinity, and therefore the sequence of the reciprocals ---that bounds above the sequence of the lengths--- goes to $0$.
\end{proof}

As a consequence, we have:
\begin{itemize}
\item[(i)] the map $k:\oi\setminus\Qbb\to I^\Nbb$ that associates to each irrational its kneading sequence is injective, continuous, and open;
\item[(ii)] if $x$ is irrational with kneading sequence $a_1,a_2,\ldots$, and $y\in\oi$, then
$$
\lim_{n\to\infty}\psi_{a_1}\cdots\psi_{a_n}(y)=x.
$$
\end{itemize}
Indeed, let $k(x)=a_1,a_2,\ldots$; by~\S(\ref{eq1}), the set of all numbers whose kneading sequence agrees with $k(x)$ up to $a_n$ is contained in the closed unimodular interval $\Gamma_n$ whose extrema are $\psi_{a_1}\cdots\psi_{a_n}(0)$ and 
$\psi_{a_1}\cdots\psi_{a_n}(1)$. Since $D$ contains at least three points, the inclusions $\Gamma_1\supset\Gamma_2\supset\cdots$ are proper, and Observation~\ref{ref1} implies~(ii) and the injectivity of $k$. We leave the continuity and openness of $k$ as an exercise for the reader; $k$ may be surjective (e.g, in the case of the ordinary continued fractions), but is never so if $I$ is finite (because then $I^\Nbb$ is compact, while $\oi\setminus\Qbb$ is not).

\section{$\ldots$ and periodicity}

On August 25th 1769, Lagrange read at the Royal Berlin Academy of Science a \emph{M\'emoire} on the resolution of algebraic equations via continued fractions~\cite{lagrange1770}. He was in Berlin since 1766, as director of the mathematical section of the same Academy, on the recommendation of his predecessor, Leonhard Euler. In this \emph{M\'emoire} he proved that the irrational $x$ has an eventually periodic continued fractions expansion iff it has degree two over the rationals. The ``only if'' implication is clear even in our general setting. If $x$ has kneading sequence $\vect a1t,\overline{\vect a{{t+1}}{{t+r}}}$ under any Gauss map $G$, then $G^tx$ and $G^{t+r}x$ have the same kneading sequence, namely $\overline{\vect a{{t+1}}{{t+r}}}$, and hence are equal. This implies that the vectors
$$
G_{a_t}\cdots G_{a_1}
\begin{pmatrix}
x\\1
\end{pmatrix},
$$
and
$$
G_{a_{t+r}}\cdots G_{a_{t+1}}G_{a_t}\cdots G_{a_1}
\begin{pmatrix}
x\\1
\end{pmatrix},
$$
are projectively equal, i.e., differ by a nonzero multiplicative constant. But then $(x\;1)^{tr}$ is an eigenvector for the integer matrix $(G_{a_t}\cdots G_{a_1})^{-1}G_{a_{t+r}}\cdots G_{a_{t+1}}G_{a_t}\cdots G_{a_1}$, and hence $x$ is quadratic.

The ``if'' direction is trickier.
The usual proof~\cite{hardywri85},
\cite{rockettszusz92}, is basically Lagrange's. One considers the minimal polynomial $c_nX^2+d_nX+e_n\in\Zbb[X]$ of $1/G^nx$, and shows that, from some $n$ on, the coefficients $c_n,d_n,e_n$ must satisfy certain inequalities. Sometimes this proof is reworded by 
writing $1/G^nx$ in reduced form
$$
\frac{1}{G^nx}=\frac{P_n+\sqrt{D}}{Q_n},
$$
and again bounding $P_n$ and $Q_n$ in terms of the common discriminant $D$ of the above polynomials, for $n$ large enough. All of this is tightly related to Gauss's theory of reduced quadratic forms;
see~\cite[\S5.7]{cohen93}.

It is not clear how to adapt the above proof to the case of arbitrary piecewise-fractional expansions. I will switch to a more geometric vein by formulating another observation.

\begin{observation}
Let $x$ be a quadratic irrational. Then there exists a matrix $H$ with integer entries having $\vc{x}{1}$ as an eigenvector, with corresponding eigenvalue $\lambda>\bar\lambda>0$ ($\bar\lambda$ being the other eigenvalue).
\end{observation}
\begin{proof}
Let $cX^2+dX+e\in\Zbb[X]$ be the minimal polynomial of $x$, and let $\bar x$ be the algebraic conjugate. Plainly
$$
\begin{pmatrix}
-d & -e \\
c & 0
\end{pmatrix}
\begin{pmatrix}
x\\ 1
\end{pmatrix}=
cx
\begin{pmatrix}
x\\ 1
\end{pmatrix}.
$$
Let $K$ be the sum of the matrix displayed above with $t$ times the $2\times2$ identity matrix. For $t$ a sufficiently large positive integer, the eigenvalues $cx+t$ and $c\bar x+t$ of $K$ are both $>0$. If $c\bar x+t < cx+t$, we take $H=K$ and we are through. Otherwise, we take $H=\abs{K}K^{-1}$ and observe that $0<\abs{K}(c\bar x+t)^{-1}=cx+t<
c\bar x+t=\abs{K}(cx+t)^{-1}$.
\end{proof}

We may now prove the Lagrange Theorem for arbitrary piecewise-fractional expansions: let $x\in\oi$ be a quadratic irrational, with kneading sequence $a_1,a_2,\ldots$ under some Gauss map $G$, and let $H,\lambda,\bar\lambda,\bar x$ be as above. Let $h:\Rbb^2\to\Rbb^2$ be the affine transformation whose matrix is $H$, with respect to the standard basis.
By~\S1(\ref{eq1}), $x$ belongs to the interval $\Gamma_n$ whose extrema are $\psi_{a_1}\cdots\psi_{a_n}(0)$ and $\psi_{a_1}\cdots\psi_{a_n}(1)$, for every $n\ge0$.
The vectors $\vc{\psi_{a_1}\cdots\psi_{a_n}(0)}{1}$ and 
$\vc{\psi_{a_1}\cdots\psi_{a_n}(1)}{1}$ are positively proportional to the columns of the matrix
$$
B_n=G_{a_1}^{-1}G_{a_2}^{-1}\cdots G_{a_n}^{-1}
\begin{pmatrix}
0 & 1 \\
1 & 1
\end{pmatrix}.
$$
Let $C_n$ be the cone spanned positively by these columns.
By \S2(ii), there exists $m$ such that $\bar x\notin\Gamma_n$, for every $n\ge m$.
Since $0<\bar\lambda<\lambda$, $\vc{x}{1}\in C_n$, and $\vc{\bar x}{1}\notin C_n\cup-C_n$, we have $h[C_n]\subseteq C_n$ for $n\ge m$. This implies that the matrix $H_n=B_n^{-1}HB_n$ (i.e., the matrix of $h$ with respect to the basis given by the columns of $B_n$) has nonnegative entries for $n\ge m$. Hence $H_m,H_{m+1},\ldots$ are nonnegative integer matrices, all conjugate to each other via matrices in $\GL_2\Zbb$. Since the determinant and the trace of a matrix are invariant under conjugation, these matrices are finite in number. Therefore $H_t=H_{t+r}$, for some $t\ge m$ and \hbox{$r>0$}. The $\lambda$-eigenspace of $H_t$ is $1$-dimensional with basis $B_t^{-1}\vc{x}{1}$, and analogously for $H_{t+r}$. Multiplying to the left by $\bigl( \begin{smallmatrix}
0&1\\ 1&1
\end{smallmatrix} \bigr)$, we see that the vectors
$$
G_{a_t}\cdots G_{a_1}
\begin{pmatrix}
x\\1
\end{pmatrix},
$$
and
$$
G_{a_{t+r}}\cdots G_{a_{t+1}}G_{a_t}\cdots G_{a_1}
\begin{pmatrix}
x\\1
\end{pmatrix},
$$
are projectively equal.
This means $G^tx=G^{t+r}x$; hence $x$ is preperiodic under the Gauss map, and its kneading sequence is eventually periodic.

\end{document}